\documentclass[a4paper,twoside,11pt,reqno]{article}

\usepackage[plain]{fullpage} 
\usepackage{amsmath}
\usepackage{amssymb}
\usepackage{latexsym}
\usepackage{amsxtra}
\usepackage{amscd}
\usepackage{theorem}
\usepackage{mathabx}
\input xy
\xyoption{all}
\entrymodifiers={+!!<0pt,\fontdimen22\textfont2>}

\usepackage{graphics}
\usepackage{epic}

\usepackage{mathrsfs, amsfonts, latexsym}

\usepackage{units}

\usepackage{color}

\theoremstyle{change}

{\theorembodyfont{\slshape}

\newtheorem{thm}{Theorem.}[section]
\newtheorem{cor}[thm]{Corollary.}
\newtheorem{lem}[thm]{Lemma.}
\newtheorem{prop}[thm]{Proposition.}

{\theorembodyfont{\rmfamily}

\newtheorem{rem}[thm]{Remark.}

}

\renewcommand{\em}{\sl}

\DeclareMathOperator*{\colim}{colim}

\parindent1em
\parskip0.3ex

\newcommand{\proof}{\noindent {\bf Proof:\ }}
\newcommand{\Endproof}{\hspace*{\fill} $\Box$ \vspace{1ex} \noindent }

\makeatletter
\renewcommand{\subsection}{\@startsection{subsection}{2}%
{\z@}{-3.25ex plus -1ex minus-.2ex}{-1em}{\bf}} \makeatother

\newcommand{\et}{{\rm\acute{e}t}}

\newcommand{\LeftEqNo}{\let\veqno\leqno}

\numberwithin{equation}{section}
\numberwithin{thm}{section}
\theoremstyle{plain}

\title{Homotopy Rational Points of Brauer-Severi Varieties}
\author{Johannes Schmidt\thanks{Supported by DFG-Forschergruppe 1920 "Symmetrie, Geometrie und Arithmetik", Heidelberg--Darmstadt}}

\date{\vspace{-5ex}}

\begin{document}

\maketitle

\begin{quotation} 
\noindent \small {\bf Abstract.}
We study homotopy rational points of Brauer-Severi varieties over fields of characteristic zero.
We are particularly interested if a Brauer-Severi variety admitting a homotopy rational point splits.
The analogue statement turns out to be true for open subvarieties of twisted hyperplane arrangements.
In the complete case, the obstruction turns out to be the pullback of the Chern class of a generator of the Picard group along the homotopy rational point.
Moreover, we will show that every Brauer-Severi variety over a $p$-adic field with $p$-primary period admits a homotopy rational point. 
\end{quotation}

\section*{Introduction}

\subsection*{Background.}
Let $X$ be a geometrically connected variety over a field $k$.
The structure morphism induces a map $\hat{\rm Et}(X) \rightarrow \hat{\rm Et}(k)$ on profinite \'{e}tale homotopy types and each $k$-rational point on $X$ induces a splitting of this map.
We call a splitting of $\hat{\rm Et}(X) \rightarrow \hat{\rm Et}(k)$ in the homotopy category of simplicial profinite sets (cf.\ \cite{Quick08}) a \textbf{homotopy rational point}.
The immediate question arising is:
\begin{center}
 (\dag)~ \parbox[t]{14cm}{\textit{Does the existence of a homotopy rational point imply the existence of a genuinely rational point?}}
\end{center}
If $X$ is of type $K(\pi,1)$, then the question for homotopy rational points is (after choosing base points) equivalent to that for splittings of the induced exact sequence of \'{e}tale fundamental groups
\begin{equation*}\LeftEqNo\tag*{$\pi_1(X/k)$:}
 \xymatrix{
  \boldsymbol{1} \ar[r] &
  \pi_1(X\otimes_kk^{\rm s},\bar{x}) \ar[r] &
  \pi_1(X,\bar{x}) \ar[r] &
  \Gamma_k \ar[r] &
  \boldsymbol{1}
 }.
\end{equation*}
In the case of hyperbolic curves over absolutely finitely generated (or $p$-adic) fields, this question is known under the name ($p$-adic) \emph{Section conjecture of anabelian geometry} -- and is still open.
In the case of genus $1$ curves over $p$-adic fields, (\dag) has a negative answer by \cite{Stix12} Prop.\ 183.

\subsection*{Aim and Results of the paper.}
As $K(\pi,1)$-spaces, smooth projective curves of positive genus are simple from a homotopical point of view.
In contrast, we want to study (\dag) in one of the simplest but nevertheless non-trivial situations from a geometric point of view:
We are interested in (\dag) in the case of \textbf{Brauer-Severi varieties} over a field $k$ of characteristic $0$, i.e.\ the case of twisted $k$-forms of projective spaces.
In this case a question equivalent to (\dag) is whether a Brauer-Severi variety $X$ admitting a homotopy rational point \textbf{splits}, i.e., is isomorphic to a projective space $\mathbb{P}^n$ as a $k$-variety.
In contrast to the $K(\pi,1)$ case, $\pi_1(X/k)$ contains no information about rational points:
Since $X\otimes_kk^{\rm s}$ is simply connected, $\pi_1(X/k)$ degenerates to an isomorphism $\pi_1^\et(X) \stackrel{\sim}{\rightarrow} \Gamma_k$.  
\\
Up to the identification of twisted linear subvarieties, Brauer-Severi varieties are classified by their Brauer class in ${\rm Br}(k)$, where the trivial class corresponds exactly to projective spaces.
Thus, a possible strategy to deal with (\dag) is, to first identify the Brauer class of $X$ using homotopy-invariant data and then use the existence of a homotopy rational point to show that this class has to be trivial. 
The simplest example of this strategy is the following reformulation of Amitsur's Theorem (see e.g.\ \cite{GilleSzamuely06} Thm.\ 5.4.1) using the homotopy invariance of Brauer groups over fields:
Let $K$ be the function field of $X$.
Then $X$ splits over $k$, if and only if the canonical map on absolute Galois groups $\Gamma_K \rightarrow \Gamma_k$ splits, i.e., Brauer-Severi varieties satisfy an analogue of the birational weak section conjecture.
Unfortunately, in higher dimensions, Brauer groups are no longer homotopy invariant in general, so we have to work more carefully to possibly find a positive answer to (\dag).

In case of a Brauer-Severi curve $X$, (\dag) is true at least for an open affine subvariety $U$ by \cite{Stix12} Prop.\ 187.
This can be proven according to the above strategy: 
By an observation of Stix, the Brauer class of $X$ in ${\rm Br}(k)$ can be found in the Hochschild-Serre spectral sequence with local coefficients $\boldsymbol{\mu}_2$ given by the covering $U\otimes_kk^{\rm s} / U$ (see the first statement of loc.\ cit.\ Prop.\ 188).
In fact, the arguments can be extended to arbitrary dimensions to give a positive answer to (\dag) in a quasi-affine case (for the exact statement, see Thm.\ \ref{thm: weak section conjecture for affine brauer-severi varieties}, below):

\bigskip\noindent {\bf Theorem A.} {\it
Let $k$ be a field of characteristic $0$ and $X/k$ a Brauer-Severi variety.
Let $U$ be an open subvariety inside the complement of a twisted hyperplane arrangement of $X$.
Assume $U$ admits a homotopy rational point.
Then $X$ splits and $U$ admits $k$-rational points.
{\it}}

\bigskip
Let us turn our attention to (\dag) in the case of a complete Brauer-Severi $X$.
Say $X$ has period $d$, i.e., the Brauer class of $X$ has order $d$ in ${\rm Br}(k)$.
The Picard group ${\rm Pic}(X)$ is the subgroup of index $d$ in ${\rm Pic}(X\otimes_kk^{\rm s}) \cong \mathbb{Z}$.
Let $[\mathcal{L}_X]$ be a generator.
We consider its image in $H^2(k,\boldsymbol{\mu}_d)$ under the mod-$d$ reduction of the first profinite Chern class map $\hat{c}_1 : {\rm Pic}(X) \rightarrow H^2(X,\hat{\mathbb{Z}}(1))$.
A homotopy rational point $s$ induces a pullback map $s^*: H^2(X,\boldsymbol{\mu}_d) \rightarrow H^2(k,\boldsymbol{\mu}_d)$.
It turns out that the pullback along $s$ of $\hat{c}_1[\mathcal{L}_X]$ generates the same subgroup of $H^2(k,\boldsymbol{\mu}_d) = {\rm Br}(k)[d]$ as the Brauer class of $X$.
Thus, we can give a positive answer to (\dag) under the according additional (and necessary) assumption on the first Chern class map $\hat{c}_1$ of $X$ (for the exact statement, see Thm.\ \ref{thm: weak section conjecture for brauer severi varieties} (ii), below):

\bigskip\noindent {\bf Theorem B.} {\it
Let $k$ be a field of characteristic $0$ and $X/k$ a Brauer-Severi variety of period $d$ admitting a homotopy rational point $s$.
Then $X$ is non-split if and only if $s^*\hat{c}_1[\mathcal{L}_X]$ is not divisible by $d$ in $H^2(k,\hat{\mathbb{Z}}(1))$.
In particular, $X$ splits if $s^*$ trivializes $\hat{c}_1$ in $H^2(k,\hat{\mathbb{Z}}(1))$.
{\it}}

\bigskip
An analogue of the section conjecture itself fails:
Over suitable $p$-adic fields, a homotopy rational point not induced by any rational point can be constructed even for projective spaces (see Rem.\ \ref{rem: projective spaces have more then just one ho rat point}, below).
In fact, the obstruction for the rationality of a homotopy rational point $s$ of $X$ is exactly the class $s^*\hat{c}_1[\mathcal{L}_X]$ in $H^2(k,\hat{\mathbb{Z}}(1))$ (see Thm.\ \ref{thm: weak section conjecture for brauer severi varieties} (iii), below):  

\bigskip\noindent {\bf Theorem C.} {\it
Let $k$ be a field of characteristic $0$ and cohomological dimension $\leq 2$ and $X/k$ a Brauer-Severi variety.
Then a homotopy rational point $s$ of $X$ is induced by a $k$-rational point if and only if $s^*$ trivializes $\hat{c}_1$ in $H^2(k,\hat{\mathbb{Z}}(1))$.
{\it}}

\bigskip
A link to question (\dag) in the cases of curves of positive genus is the following:
Let $C/k$ be a smooth projective curve of positive genus.
For a class $[A]$ in the relative Brauer group
\begin{equation*}
 {\rm Br}(C/k) := {\rm ker}(
 \xymatrix{
  {\rm Br}(k) \ar[r] &
  {\rm Br}(C)
 }),
\end{equation*}
there is a non-trivial map $f: C \rightarrow X_A$ to the Brauer-Severi variety $X_A$ corresponding to the central simple algebra $A$.
Thus, any homotopy rational point $s$ of $C$ induces a homotopy rational point $f_*s$ of $X_A$ and one could ask whether $X_A$ splits, i.e., whether $[A]$ is trivial.
Applying this construction to the genus $1$ obstruction for the curve case of (\dag) reveals, that every (non-split) Brauer-Severi variety over a $p$-adic field with $p$-primary period admits a homotopy rational point (see Prop.\ \ref{prop: counterexample}, below).
Thus, not every homotopy rational point of a Brauer-Severi variety satisfies the extra assumption on the first Chern class map of Thm.\ A.
In particular, not every homotopy rational point comes from a rational point. 
We hope that these are essentially the only obstructions of question (\dag) for Brauer-Severi varieties over $p$-adic fields:
\textit{A Brauer-Severi variety over a $p$-adic field should admit a homotopy rational point if and only if it has a $p$-primary period.}

\subsection*{Notation.}
In the following, $k$ always denotes a field of characteristic $0$.
If $K$ is any field, we denote a fixed separable closure by $K^{\rm s}$ and the corresponding absolute Galois group by $\Gamma_K$.
For $X$ a $K$-variety, denote its base-change to $K^{\rm s}$ by $X^{\rm s}$.
Denote by $\hat{\mathcal{S}}_{(*)}$ the (pointed) category of profinite simplicial sets together with the model structures of \cite{Quick08}.
For $X$ a scheme together with a geometric point $\bar{x}$, $\pi_1(X,\bar{x})$ denotes its profinite \'{e}tale fundamental group.
Denote by $\hat{\rm Et}(X)$ its profinite \'{e}tale homotopy type in $\hat{\mathcal{S}}_{(*)}$.
If no confusion could arise, we will just write $X$ for $\hat{\rm Et}(X)$.
Finally, we will always use continuous \'{e}tale cohomology in the sense of \cite{Jannsen88}.

\subsection*{Acknowledgements.}
This paper is mostly based on my dissertation, and I would like to thank my advisor, Jakob Stix, for his encouragement and guidance.
Thanks also to Armin Holschbach, Gereon Quick and Alexander Schmidt for helpful discussions and suggestions.

\section{Homotopy rational Points}

Let $X$ be a $K$-variety.
The profinite \'{e}tale homotopy type of the spectrum of $K$ is a $K(\pi,1)$ with fundamental group $\Gamma_K$, i.e., it is weakly equivalent to the profinite classification space $B\Gamma_K$ in $\hat{\mathcal{S}}_{(*)}$ (pointed by our choice of $K^{\rm s}/K$).
The structural map of $X/K$ gives a canonical map $X \rightarrow B\Gamma_K$ of profinite \'{e}tale homotopy types, i.e., we treat the homotopy type of $X/K$ as an object in $\hat{\mathcal{S}}_{(*)}\downarrow B\Gamma_K$.
The homotopy type of a $K$-rational point of $X$ is a splitting of the map $X\rightarrow B\Gamma_K$, i.e., a map $B\Gamma_K \rightarrow X$ in $\hat{\mathcal{S}}_{(*)}\downarrow B\Gamma_K$.
A {\bf homotopy rational point} of $X$ over $K$ (in $\mathcal{H}(\hat{\mathcal{S}}\downarrow B\Gamma_K)$) is a homotopy class in
\begin{equation*}
[B\Gamma_K,X]_{\hat{\mathcal{S}}\downarrow B\Gamma_K}.
\end{equation*}
We call a homotopy rational point {\bf rational}, if it is the homotopy class of an actual $K$-rational point.

Similar to the case of simplicial sets, for a profinite space $\mathfrak{X}$ and a profinite group $\pi$, homotopy classes $[\mathfrak{X},B\pi]_{\hat{\mathcal{S}}}$ correspond to functors of profinite categories modulo natural transformations between the fundamental groupoid $\Pi(\mathfrak{X})$ and $\pi$.
For $X$ a $K(\pi,1)$-space, i.e., weakly equivalent to $B\pi_1(X)$, it follows that sections of the fundamental group sequence $\pi_1(X/K)$ 
modulo conjugation correspond to homotopy rational points of $X$ over $K$ (cf.\ \cite{Quick15} 2.4).
By \cite{Stix02} Prop.\ A.4.1, this in particular is the case for $X$ any smooth curve over $K$ except for Brauer-Severi curves.

By \cite{Quick10} Cor.\ 2.11, $\mathcal{H}(\hat{\mathcal{S}}\downarrow B\Gamma_K)$ is Quillen equivalent to $\mathcal{H}(\hat{\mathcal{S}}_{\Gamma_K})$.
It follows that homotopy rational points of $X$ over $K$ can be treated as homotopy classes in $[E\Gamma_K,X\times_{B\Gamma_K}E\Gamma_K]_{\hat{\mathcal{S}}_{\Gamma_K}}$, i.e., as connected components of Quick's homotopy fixed point space
\begin{equation*}
 (X\times_{B\Gamma_K}E\Gamma_K)^{h\Gamma_K} = S_{\Gamma_K}(E\Gamma_K,\mathfrak{X})
\end{equation*}
($\mathfrak{X}$ a fibrant replacement of the base extension $X\times_{B\Gamma_K}E\Gamma_K$ in $\hat{\mathcal{S}}_{\Gamma_K}$) defined and studied in \cite{Quick10} (also cf.\ \cite{Quick15} 2.3).

In general, homotopy fixed points under a profinite group $\Gamma$ are difficult to describe.
For Eilenberg-MacLane spaces there is the following description (see \cite{Quick10} Prop.\ 2.15):
For $\Lambda$ a profinite $\Gamma$-module, the homotopy groups of $K(\Lambda,n)^{h\Gamma}$ are given by
\begin{equation*}
 \pi_q(K(\Lambda,n)^{h\Gamma}) = H^{n-q}(\Gamma,\Lambda) {\rm ~(}=0{\rm ~for~}n-q<0 {\rm )}.
\end{equation*}
For $\mathfrak{X}$ a general profinite $\Gamma$-space, there is a Bousfield-Kan type descent spectral sequence
\begin{equation}\label{eq: descent spectral sequence}
 E_2^{p,q} = H^p(\Gamma,\pi_{-q}(\mathfrak{X})) \Rightarrow \pi_{-(p+q)}(\mathfrak{X}^{h\Gamma})
\end{equation}
(with differentials in the usual ``cohomological'' directions) by \cite{Quick10} Thm.\ 2.16.
Applying Bousfield and Kan's connectivity lemma \cite{BousfieldKan} Ch.\ IX 5.1 to (\ref{eq: descent spectral sequence}) gives us a little grasp on homotopy fixed point sets.
An example for this strategy is the following lemma:

\begin{lem}\label{lem: injection of hofixed points}
Let $\Gamma$ be a profinite group of cohomological dimension $\leq n$ and $f:\mathfrak{X} \rightarrow \mathfrak{Y}$ an $(n+1)$-equivalence in $\hat{\mathcal{S}}_\Gamma$ (i.e., $\pi_q(f)$ is an isomorphism for all $q \leq n$ and an epimorphism for $q = n+1$).
Then $f$ induces an injection
\begin{equation*}
 [E\Gamma,\mathfrak{X}]_{\hat{\mathcal{S}}_\Gamma} =
 \xymatrix{
  \pi_0(\mathfrak{X}^{h\Gamma}) \ar[r] &
  \pi_0(\mathfrak{Y}^{h\Gamma})
 }
 = [E\Gamma,\mathfrak{Y}]_{\hat{\mathcal{S}}_\Gamma}.
\end{equation*}
\end{lem}

\proof
We may assume that $\mathfrak{Y}$ is fibrant and $f$ is a fibration in $\hat{\mathcal{S}}_\Gamma$.
Further, we may assume that $\mathfrak{X}^{h\Gamma}$ is non-empty.
Say, $s:E\Gamma \rightarrow \mathfrak{X}$ is a model of a homotopy fixed point and let $r$ be $f \circ s$.
The fibre $\mathfrak{F}_s:= \mathfrak{X} \times_{\mathfrak{Y}} E\Gamma$ comes equipped with a fibration into $E\Gamma$, hence is fibrant in $\hat{\mathcal{S}}_\Gamma$, too.
Taking limits (i.e., forgetting the topology in \cite{Quick10}) resp.\ simplicial mapping spaces $S_\Gamma(E\Gamma,-)$ of $\hat{\mathcal{S}}_\Gamma$ gives us a homotopy fibre sequence
\begin{equation*}
 \xymatrix{
  \lim\mathfrak{F}_s \ar[r] &
  \lim\mathfrak{X} \ar[r] &
  \lim\mathfrak{Y}
 }
\end{equation*}
resp.\
\begin{equation*}
 \xymatrix{
  \mathfrak{F}_s^{h\Gamma} \ar[r] &
  \mathfrak{X}^{h\Gamma} \ar[r] &
  \mathfrak{Y}^{h\Gamma}
 }
\end{equation*}
in $\underline{\rm SSets}_\bullet$ (pointed by the neutral element in $\Gamma$).
By \cite{Quick13} Lem. 2.9, the limit of $f$ is an $n$-equivalence of simplicial sets. 
So, again by loc.\ cit., the first fibre sequence implies the $n$-connectedness of $\mathfrak{F}_s$.
Using the second homotopy fibre sequence, we get that the map of pointed sets $(\pi_0(\mathfrak{X}^{h\Gamma}),s) \rightarrow (\pi_0(\mathfrak{Y}^{h\Gamma}),r)$ has kernel $\pi_0(\mathfrak{F}_s^{h\Gamma})$.
Bousfield and Kan's connectivity lemma applied to the descent spectral sequence (\ref{eq: descent spectral sequence}) for $\mathfrak{F}_s$ implies that this kernel is trivial, since $\mathfrak{F}_s$ is $n$-connected and $\Gamma$ has cohomological dimension $\leq n$.
Varying over all the homotopy fixed points of $\mathfrak{X}$, we get the result.
\Endproof

Returning to homotopy rational points, let us mention the following lemma about homotopy classes of $K$-rational points:

\begin{lem}\label{lem: ho classes of rational points}
Let $Z$ be a non-empty irreducible $K$-variety.
\begin{enumerate}
 \item If $K$ has cohomological dimension $\leq n$ and $Z$ is geometrically $n$-connected, then the canonical map
 \begin{equation*}
  \xymatrix{
   Z(K) \ar[r] &
   [B\Gamma_K,Z]_{\hat{\mathcal{S}}\downarrow B\Gamma_K}
  }
 \end{equation*}
 is trivial.  
 \item If $K$ has characteristic $0$ and $Z$ is $\mathbb{A}^1$-chain connected (i.e., for any two $K$-rational points $z^\prime,z^{\prime\prime}$, there are finitely many $K$-morphisms $u_i: \mathbb{A}_K^1 \rightarrow Z$, $1\leq i \leq n$ with $u_1(0) = z^\prime$, $u_n(1) = z^{\prime\prime}$ and $u_i(1) = u_{i+1}(0)$), 
 then the canonical map
 \begin{equation*}
  \xymatrix{
   Z(K) \ar[r] &
   [B\Gamma_K,Z]_{\hat{\mathcal{S}}\downarrow B\Gamma_K}
  }
 \end{equation*}
 is trivial.  
\end{enumerate}
\end{lem}

\proof
To prove (i), it is enough to show that $\pi_0((Z\times_{B\Gamma_K}E\Gamma_K)^{h\Gamma_K})$ is trivial.
But this follows from Bousfield and Kan's connectivity lemma applied to the descent spectral sequence (\ref{eq: descent spectral sequence}) for $Z\times_{B\Gamma_K}E\Gamma_K$:
${\rm cd}(K) \leq n$ and $Z$ is geometrically $n$-connected by assumption.
Statement (ii) holds, since $\mathbb{A}_K^1$ is contractible (to $B\Gamma_K$) in characteristic $0$.
\Endproof

\section{Homotopy invariants of Brauer-Severi varieties}

\subsection*{Brauer-Severi varieties.}
Recall, that a Brauer-Severi variety over $k$ is a $k$-variety whose base extension to an algebraic closure $k^{\rm s}/k$ of $k$ is $k^{\rm s}$-isomorphic to the projective space $\mathbb{P}^n$ for a suitable $n$.
It is a basic fact that a Brauer-Severi variety has a $k$-rational point if and only if it is isomorphic over $k$ to a projective space (see e.g.\ \cite{GilleSzamuely06} Thm. 5.1.3).
We say that $X$ \textbf{splits} in this case.
Let us first discuss the fundamental algebraic topology of Brauer-Severi varieties:

\begin{rem}\label{rem: non equivariant etale algebraic topology of projective spaces}
Let $X$ be a Brauer-Severi variety over $k$ (as always: of characteristic $0$).
The higher homotopy groups of $X$ can be computed over $k^{\rm s}$.
By \cite{ArtinMazur} Cor.\ 12.12, we may replace $k^{\rm s}$ by $\mathbb{C}$.
Using the generalized Riemann Existence Theorem (loc.\ cit.\ Cor.\ 12.10), we find (without $\Gamma_k$-action):
\begin{equation*}
 \pi_q(X) \cong
 \begin{cases}
  \Gamma_k &
  {\rm if}~ q=1
 \\
  \hat{\mathbb{Z}} &
  {\rm if}~ q=2
 \\
  \pi_q(S^{2n+1}) &
  {\rm if}~ q\neq1,2.
 \end{cases}
\end{equation*}
\end{rem}

Unfortunately, this description does not include the $\Gamma_k$-action on the higher homotopy groups.
Since this action is just the canonical Galois-action given by functoriality on $\pi_q(X \otimes_k k^{\rm s})$, we can reduce an equivariant description to the case of a projective space:

\begin{lem}\label{lem: Galois representations on homotopy invariants}
Let $X$ and $Y$ be two Brauer-Severi varieties over $k$ of the same dimension.
Then the induced $\Gamma_k$-objects in the homotopy categorie $\mathcal{H}(\hat{\mathcal{S}})$ are canonically isomorphic.
In particular, all Galois representations on geometric homotopy invariants of $X$ and $Y$ (e.g., $\ell$-adic cohomology, higher homotopy groups, \dots) are canonically isomorphic. 
\end{lem}

The proof will turn out to be quite easy after the following two remarks: 

\begin{rem}\label{rem: Eilenberg-Zilber}
Let $X$ be a geometrically unibranch proper and simply connected $k^{\rm s}$-variety and let $Y$ be any normal $k^{\rm s}$-variety.
Then the canonical map
\begin{equation*}
 \xymatrix{
 \hat{\rm Et}(X \times Y) \ar[r] &
 \hat{\rm Et}(X) \times \hat{\rm Et}(Y)
 }
\end{equation*}
is a weak equivalence in $\hat{\mathcal{S}}_*$:
It suffices to check that the canonical morphisms from $\pi_q$ of the source resp.\ target into $\pi_q(X) \times \pi_q(Y)$ are isomorphisms.
For the target, this follows from the long exact sequence of homotopy groups of a fibration of simplicial sets.
For the source, we use that all occurring homotopy groups are profinite together with the analogue sequence of \cite{Friedlander73} Cor.\ 4.8 applied to the projection $X \times Y \rightarrow Y$.
\end{rem}

\begin{rem}\label{rem: trivial group action}
Let $G$ be a connected, quasi-projective group scheme over $k^{\rm s}$ and $\mu: X \times G \rightarrow X$ a proper, simply connected and geometrically unibranch $G$-space in $\underline{\rm Var}_{k^{\rm s}}$.
Then the induced $G(k^{\rm s})$-action on $X$ is trivial in $\mathcal{H}(\hat{\mathcal{S}})$:
The $G(k^{\rm s})$-action on $X$ is given by
\begin{equation*}
X=
 \xymatrix{
  X\times {\rm Spec}(k^{\rm s}) \ar[r]^-{{\rm id} \times g} \ar@/_1pc/[rr]_-{(-).g} &
  X\times G \ar[r]^-\mu &
  X
 }
\end{equation*}
for $g \in G(k^{\rm s})$.
Thus, the claim follows from Rem.\ \ref{rem: Eilenberg-Zilber} and Lem.\ \ref{lem: ho classes of rational points} (i).
\end{rem}

\noindent \textbf{Proof of Lem.\ \ref{lem: Galois representations on homotopy invariants}:}
We may assume that $Y$ splits, i.e., is isomorphic to $\mathbb{P}_k^n$.
Two trivializations $f,g: X^{\rm s} \cong \mathbb{P}_{k^{\rm s}}^n$
over $k^{\rm s}$ differ by the element $h = g \circ f^{-1}$ in ${\rm PGL}_{n+1}(k^{\rm s})$.
By Rem.\ \ref{rem: trivial group action} they agree in $\mathcal{H}(\hat{\mathcal{S}})$ and thus define a canonical isomorphism in the respective homotopy category.
But for any such trivialization $f$ and $\gamma \in \Gamma_k$, the translate $\gamma \circ f \circ \gamma^{-1}$ is another such trivialization, so  our canonical isomorphism is even $\Gamma_k$-equivariant.
\Endproof

\begin{rem}\label{rem: continuous vs abstract action}
The higher homotopy or reduced homology groups of a Brauer-Severi variety are profinite continuous $\Gamma_k$-modules but Lem.\ \ref{lem: Galois representations on homotopy invariants} only compares the abstract $\Gamma_k$-actions on the induced profinite abelian groups.
To deal with this subtlety, note that the category of profinite continuous $\Gamma_k$-modules is a full subcategory of the category of profinite abelian groups with an abstract $\Gamma_k$-action:
This follows using Pontryagin duality, since the category of abstract $\Gamma_k$-actions on torsion abelian groups contains the category of torsion discrete continuous $\Gamma_k$-modules as a full subcategory. 
\end{rem}

\begin{rem}\label{rem: elements of equivariant etale algebraic topology of projective spaces}
The geometric cohomology with coefficients $\mathbb{Z}/n\mathbb{Z}(m)$ of a projective space is well known.
Using the universal coefficients theorem and Lem.\ \ref{lem: Galois representations on homotopy invariants}, for $X$ an $n$-dimensional Brauer-Severi variety over $k$ we find:   
\begin{equation*}
 \tilde{H}_q(X^{\rm s}) =
 \begin{cases}
  0 &
  {\rm if}~ q {\rm~is~odd,~}>2n{\rm~or~}=0
 \\
  \hat{\mathbb{Z}}(\frac{q}{2}) &
  {\rm else}
 \end{cases}
\end{equation*}
In particular, $\pi_2(X)$ is isomorphic to $\hat{\mathbb{Z}}(1)$ as a $\Gamma_k$-module by Hurewicz.
\end{rem}

Finally, let us discuss the second absolute cohomology of a Brauer-Severi variety:

\begin{lem}\label{lem: second cohomology of Brauer-Severi varieties}
Let $X$ be a Brauer-Severi variety over $k$ admitting a homotopy rational point $s$ and let $p:\ X^{\rm s} \rightarrow X$ be the canonical projection.
Let $\Lambda$ be a profinite continuous $\Gamma_k$-module. 
Then we get a direct sum decomposition
\begin{equation*}
 H^2(X,\Lambda) \cong
 H^2(B\Gamma_k,\Lambda) \oplus
 H^2(X^{\rm s},\Lambda),
\end{equation*}
where the projections onto the two summands correspond to $s^*$ resp.\ $p^*$.
\end{lem}

\proof
Consider the Hochschild-Serre spectral sequence ${\rm HS}_*^{\bullet,\bullet}(X,\Lambda)$ in \cite{Jannsen88} Cor.\ 3.4:
We claim that the differential $\partial_3^{0,2}: {\rm HS}_3^{0,2}(X,\Lambda) \rightarrow  {\rm HS}_3^{3,0}(X,\Lambda)$
is trivial:
The canonical map $H^3(\Gamma_k,\Lambda) \rightarrow H^3(X,\Lambda)$ factors as
\begin{equation*}
 \xymatrix{
  H^3(\Gamma_k,\Lambda) \ar[r] \ar@{=}[d] &
  H^3(X,\Lambda)\phantom{.}
 \\
  {\rm HS}_2^{3,0}(X,\Lambda) \ar@{->>}[r] &
  {\rm HS}_\infty^{3,0}(X,\Lambda). \ar@{^(->}[u]
 }
\end{equation*}
By \cite{Quick08} Prop. 3.1., $s^*$ gives a retraction of this map.
It follows that the lower horizontal arrow is an isomorphism which forces the image of $\partial_3^{0,2}$ to be trivial (for general $k$-varieties, this argument would give $\partial_2^{1,1} = 0$, as well).
Thus, we get the split exact sequence
\begin{equation*}
 \xymatrix{
  0 \ar[r] &
  H^2(\Gamma_k,\Lambda) \ar[r] &
  H^2(X,\Lambda) \ar[r]^-{p^*} &
  H^2(X^{\rm s},\Lambda) \ar[r] &
  0
 },
\end{equation*}
where a retraction is given by $s$.
\Endproof

Using the following construction, we can make the isomorphism $\pi_2(X) = \hat{\mathbb{Z}}(1)$ more explicit: 

\begin{rem}\label{rem: fake chern class}
As a $K^{\Gamma_k}(\hat{\mathbb{Z}}(1),2)$-space (the Borel-construction of the equivariant Eilenberg-MacLane space $K(\hat{\mathbb{Z}}(1),2)$), $\mathbb{P}_k^\infty$ represents $H^2(-,\hat{\mathbb{Z}}(1))$ in $\mathcal{H}(\hat{\mathcal{S}}\downarrow B\Gamma_k)$.
The canonical embedding $\mathbb{P}_k^n \hookrightarrow \mathbb{P}_k^\infty$ corresponds to the Chern-class $\hat{c}_1[\mathcal{O}(1)]$ in $H^2(\mathbb{P}_k^n,\hat{\mathbb{Z}}(1))$.
For $X$ a Brauer-Severi variety and $s$ a homotopy rational point, there is a similar map
\begin{equation*}
 \xymatrix{
  \alpha_s: X \ar[r] &
  K^{\Gamma_k}(\hat{\mathbb{Z}}(1),2)
 }
\end{equation*}
in $\mathcal{H}(\hat{\mathcal{S}}\downarrow B\Gamma_k)$, too:
By Lem.\ \ref{lem: second cohomology of Brauer-Severi varieties} there is a unique homotopy class $\alpha_s$ corresponding to the unique cohomology class (also denoted by $\alpha_s$) with $p^*\alpha_s$ the Chern-class of $\mathcal{O}_{X^{\rm s}}(1)$ and trivial $s^*\alpha_s$.
By construction, $\alpha_x = \hat{c}_1[\mathcal{O}_X(1)]$ for $x$ a $k$-rational point of $X$. 
\end{rem}

\begin{rem}\label{rem: properties of the fake chern class}
It is immediate that $\alpha_s$ induces isomorphisms on $\pi_0$ and $\pi_1$.
By construction, ${\rm res}_{\boldsymbol{1}}^{\Gamma_k}(\alpha_s \times_{B\Gamma_k}E\Gamma_k)$ is isomorphic in $\mathcal{H}(\hat{\mathcal{S}})$ to the canonical embedding $\mathbb{P}_{k^{\rm s}}^n \hookrightarrow \mathbb{P}_{k^{\rm s}}^\infty$, which induces an isomorphism on $\pi_2$ (and on cohomology with local coefficients in degrees $\leq 2n +1$), as well.
Thus, $\pi_2(\alpha_s)$ gives back the canonical isomorphism $\pi_2(X) = \hat{\mathbb{Z}}(1)$.
\end{rem}

\subsection*{The set of homotopy rational points of a Brauer-Severi variety.}
Homotopy fixed points of Eilenberg-MacLane spaces are given by Galois-cohomology (see \cite{Quick10} Prop.\ 2.15), so the set of homotopy rational points of $K^{\Gamma_k}(\hat{\mathbb{Z}}(1),2)$ is given by $H^2(\Gamma_k,\hat{\mathbb{Z}}(1))$.
The ``fake Chern-class'' $\alpha_s$ of Rem.\ \ref{rem: fake chern class} induces an injection between the respective sets of homotopy rational points:

\begin{prop}\label{prop: injection of horat points}
Let $X$ be a Brauer-Severi variety over $k$ and suppose that $k$ has cohomological dimension $\leq 2$ (e.g.\ a $p$-adic or totally imaginary number field).
Then $\alpha_s$ induces an injection
\begin{equation*}
 \xymatrix{
  [B\Gamma_k,X]_{\hat{\mathcal{S}}\downarrow B\Gamma_k} \ar@{^(->}[r] &
  [B\Gamma_k,K^{\Gamma_k}(\hat{\mathbb{Z}}(1),2)]_{\hat{\mathcal{S}}\downarrow B\Gamma_k}
 } = H^2(\Gamma_k,\hat{\mathbb{Z}}(1)).
\end{equation*}
\end{prop}

\proof
By Rem.\ \ref{rem: properties of the fake chern class}, $\alpha_s$ induces isomorphisms on homotopy groups in degree $\leq 2$ and an epimorphism in degree $3$.
Using the Quillen equivalence \cite{Quick10} Cor. 2.11, Prop.\ \ref{prop: injection of horat points} follows from Lem.\ \ref{lem: injection of hofixed points}.
\Endproof

\section{An analogue for the weak section conjecture: complements of twisted hyperplane arrangements.}

Let $X$ be a Brauer-Severi variety.
Fix a trivialization $X^{\rm s} \cong \mathbb{P}^n\otimes_kk^{\rm s}$ over $k^{\rm s}$.
Choose hyperplanes $H_i$ of $\mathbb{P}^n\otimes_k\bar{k}$ and denote by $D^{\rm s} \hookrightarrow X^{\rm s}$ the union of the $\Gamma$-orbits of the corresponding closed subschemes of $X^{\rm s}$ (together with the reduced structure), i.e., $D^{\rm s}$ corresponds to an arrangement of hyperplanes of $\mathbb{P}^n\otimes_k\bar{k}$.
Now $\bar{D}$ descents to a hypersurface $D \hookrightarrow X$.
We call such a hypersurface a \textbf{twisted hyperplane arrangement}.

In this section, we will give a positive answer to (\dag) for open subvarieties of $X \smallsetminus D$, i.e., we will prove an analogue of the weak section conjecture for open subvarieties in complements of twisted hyperplane arrangements:

\begin{thm}\label{thm: weak section conjecture for affine brauer-severi varieties}
Let $k$ be a field of characteristic $0$ and $X/k$ a Brauer-Severi variety.
Let $U$ be an open subvariety inside the complement of a twisted hyperplane arrangement $D$ of $X$.
Assume $U$ admits a homotopy rational point $s \in [B\Gamma_k,U]_{\hat{\mathcal{S}}\downarrow B\Gamma_k}$.
Then $X$ splits and $U$ admits $k$-rational points.
\end{thm}

\begin{rem}\label{rem: weak section conjecture for affine genus 0 curves}
The one dimensional case of Thm.\ \ref{thm: weak section conjecture for affine brauer-severi varieties} was first proven by Stix (see \cite{Stix12} Prop.\ 187).
The proof of Thm.\ \ref{thm: weak section conjecture for affine brauer-severi varieties} essentially is just an extension of \cite{Stix12} Prop.\ 188 (2) to higher dimensions.
\end{rem}

To prove Thm.\ \ref{thm: weak section conjecture for affine brauer-severi varieties}, let us first make the following reductions:

\begin{rem}\label{rem: oBdA U as bis as possible}
It is enough to show that $X$ splits over $k$:
Indeed, the $k$-rational points are dense in $\mathbb{P}^n$.
In particular, we may assume that $U$ is the whole complement $X \smallsetminus D$.
\end{rem}

\begin{rem}\label{rem: amitsur theorem}
The generic point of a smooth connected $k$-variety induces a monomorphism on Brauer-groups (see e.g.\ \cite{Milne} Chapt.\ III Ex.\ 2.22).
In particular, the relative Brauer groups ${\rm Br}(X/k)$ and ${\rm Br}(U/k)$ agree and are generated by the Brauer class $[X]$ of $X$ by Amitsur's Theorem (see e.g.\ \cite{GilleSzamuely06} Thm.\ 5.4.1).
Thus, it suffices to show that ${\rm Br}(U/k)$ is trivial.
\end{rem}

The key observation for the proof of Thm.\ \ref{thm: weak section conjecture for affine brauer-severi varieties} is the following extension of \cite{Stix12} Prop.\ 188 (2) to higher dimensions:

\begin{lem}\label{lem: brauer class and hochschild serre}
Let $D$ be a twisted hyperplane arrangement of $X$ with open complement $U$.
Let $m$ be a multiple of the period of $X$ and the degree of $D^{\rm s}$.
Consider the Hochschild-Serre spectral sequence with coefficients $\boldsymbol{\mu}_m$ for the covering $U^{\rm s}/U$.
Then ${\rm Br}(U/k)$ is the image of the differential
\begin{equation*}
 \partial_2^{0,1}:
 \xymatrix{
  H^0(\Gamma,H^1(U^{\rm s},\mu_m)) \ar[r] & H^2(\Gamma,H^0(U^{\rm s},\mu_m))
 }.
\end{equation*}
\end{lem}

\proof
Denote by $F^\bullet$ the descending filtration on $H^2(U,\boldsymbol{\mu}_m)$ resp.\ $H^2(U,\mathbb{G}_{\rm m})$ of the resp.\ Hochschild-Serre spectral sequence and denote by $p$ be the canonical map $U \rightarrow {\rm Spec}(k)$.
From the $5$-term sequence we get the exact sequence
\begin{equation*}
 \xymatrix{
  H^0(\Gamma,H^1(U^{\rm s},\boldsymbol{\mu}_m)) \ar[r]^-{\partial_2^{0,1}} &
  H^2(\Gamma,\boldsymbol{\mu}_m) \ar[r] \ar[rd]_-{p^*} &
  F^1H^2(U,\boldsymbol{\mu}_m) \ar[r] \ar@{^(->}[d] &
  H^1(\Gamma,H^1(U^{\rm s},\boldsymbol{\mu}_m))
 \\
   & &
   H^2(U,\boldsymbol{\mu}_m) &
 }
\end{equation*}
In particular, $\partial_2^{0,1}$ has image ${\rm ker}(p^*)$ and the kernel of $F^1H^2(U,\boldsymbol{\mu}_m) \rightarrow H^1(\Gamma,H^1(U^{\rm s},\boldsymbol{\mu}_m))$ is the image of $p^*$.
Further, from the Kummer sequence we get the exact sequence:
\begin{equation*}
 \xymatrix{
  {\rm Pic}(U) \ar[r]^-{\cdot m} &
  {\rm Pic}(U) \ar[r]^-{\delta} &
  H^2(U,\mu_m) \ar[r] &
  {\rm Br}(U)[m] \ar[r] &
  0
 \\
   & &
  H^2(\Gamma,\boldsymbol{\mu}_m) \ar[r]^\cong \ar[u]^{p^*} &
  {\rm Br}(k)[m] \ar[u]
 }
\end{equation*}
Since the period of $X$ divides $m$, ${\rm Br}(U/k)$ lies inside ${\rm Br}(k)[m]$.
It follows that $\partial_2^{0,1}$ has image inside ${\rm Br}(U/k)$ and equality holds if and only if ${\rm im}(p^*)\cap{\rm im}(\delta) = 0$.

As $U^{\rm s}$ is isomorphic to the complement of a nonempty arrangement of hyperplanes in $\mathbb{P}_{k^{\rm s}}^n$, it is isomorphic to the spectrum of a unique factorization domain and thus, ${\rm Pic}(U^{\rm s})$ is trivial.
From Hochschild-Serre, we get ${\rm Pic}(U) = F^1H^1(U,\mathbb{G}_{\rm m}) = H^1(\Gamma_k,\mathcal{O}^\times(U^{\rm s}))$, i.e., $\delta$ factors through $F^1H^2(U,\boldsymbol{\mu}_m)$, too.
Summing up, ${\rm im}(p^*)$ and ${\rm im}(\delta)$ intersect trivially, if and only if the composition
\begin{equation*}
 \xymatrix{
   H^1(\Gamma_k,\mathcal{O}^\times(U^{\rm s})) \ar[r] \ar@/_1pc/[rr] &
   F^1H^2(U,\boldsymbol{\mu}_m) \ar[r] &
   H^1(\Gamma_k,H^1(U^{\rm s},\boldsymbol{\mu}_m))
 }
\end{equation*}
is injective.
Analyzing the Kummer sequence and relevant Hochschild-Serre spectral sequences, we see that this map is induced on $H^1(\Gamma_k,-)$ by $\delta^{\rm s}: \mathcal{O}^\times(U^{\rm s}) \rightarrow H^1(U^{\rm s},\boldsymbol{\mu}_m)$ of the Kummer sequence over $U^{\rm s}$.
Since ${\rm Pic}(U^{\rm s})$ is trivial, this latter map can be identified with the mod-$m$ quotient of $\mathcal{O}^\times(U^{\rm s})$.

Denote by ${\rm Div}_{X^{\rm s},D^{\rm s}}$ resp.\ ${\rm Div}_{X^{\rm s},D^{\rm s}}^0$ the group of divisors of $X^{\rm s}$ supported in $D^{\rm s}$ resp.\ those of degree $0$.
Since $U^{\rm s}$ is isomorphic to the spectrum of a unique factorization domain, we get a commutative diagram with exact rows:
\begin{equation*}
 \xymatrix{
  0 \ar[r] &
  k^{s,\times} \ar[r] \ar@{->>}[d]^-{(-)^m} &
  \mathcal{O}^\times(U^{\rm s}) \ar[r] \ar[d]^-{(-)^m} &
  {\rm Div}_{X^{\rm s},D^{\rm s}}^0 \ar[r] \ar[d]^-{\cdot m} &
  0
 \\
  0 \ar[r] &
  k^{s,\times} \ar[r] &
  \mathcal{O}^\times(U^{\rm s}) \ar[r] &
  {\rm Div}_{X^{\rm s},D^{\rm s}}^0 \ar[r] &
  0
 }
\end{equation*}
In particular, the mod-$m$ quotient of $\mathcal{O}^\times(U^{\rm s})$ factors through the mod-$m$ quotient of ${\rm Div}_{X^{\rm s},D^{\rm s}}^0$.
In this factorization, $\mathcal{O}^\times(U^{\rm s}) \rightarrow {\rm Div}_{X^{\rm s},D^{\rm s}}^0$ induces an injection on $H^1(\Gamma_k,-)$ by Hilbert $90$.
Thus, we are done if we can show that multiplication by $m$ on ${\rm Div}_{X^{\rm s},D^{\rm s}}^0$ induce the trivial map on $H^1(\Gamma_k,-)$.
For this, we argue as follows:
From the degree map, we get the long exact sequence on Galois cohomology:
\begin{equation*}
 \xymatrix{
  \dots \ar[r] &
  {\rm Div}_{X^{\rm s},D^{\rm s}}^{\Gamma_k} \ar[r]^-{\rm deg} &
  \mathbb{Z} \ar[r] &
  H^1(\Gamma_k,{\rm Div}_{X^{\rm s},D^{\rm s}}^0) \ar[r] &
  H^1(\Gamma_k,{\rm Div}_{X^{\rm s},D^{\rm s}}) \ar[r] &
  \dots
 }.
\end{equation*}
Now ${\rm Div}_{X^{\rm s},D^{\rm s}} = \mathbb{Z}[D^{s,(0)}]={\rm ind}_{\Gamma_k}^\Delta(\mathbb{Z})$ for a suitable open subgroup $\Delta \leq \Gamma_k$.
By Shapiro's Lemma, $H^1(\Gamma_k,{\rm Div}_{X^{\rm s},D^{\rm s}})$ is trivial so $H^1(\Gamma_k,{\rm Div}_{X^{\rm s},D^{\rm s}}^0)$ is isomorphic to $\mathbb{Z} / {\rm deg}({\rm Div}_{X,D})$.
By assumption ${\rm deg}(D^{\rm s})$ divides $m$, so $\mathbb{Z} / {\rm deg}({\rm Div}_{X,D})$ is indeed killed by $m$, finishing the proof.
\Endproof

\noindent \textbf{Proof of Thm.\ \ref{thm: weak section conjecture for affine brauer-severi varieties}:}
The section $s$ induces a morphism between Hochschild-Serre spectral sequences with coefficients $\boldsymbol{\mu}_m$ for the coverings $U^{\rm s}/U$ and $E\Gamma_k / B\Gamma_k$.
In particular, we get a commutative diagram
\begin{equation*}
 \xymatrix{
  H^0(\Gamma_k,H^1(U^{\rm s},\boldsymbol{\mu}_m)) \ar[r]^-{s^*} \ar[d]^-{\partial_2^{0,1}} &
  H^0(\Gamma_k,H^1(E\Gamma_k,\boldsymbol{\mu}_m)) \ar[d]^-{\partial_2^{0,1}} \phantom{.}
 \\
  H^2(\Gamma_k,H^0(U^{\rm s},\boldsymbol{\mu}_m)) \ar[r]^-{s^*} &
  H^2(\Gamma_k,H^0(E\Gamma_k,\boldsymbol{\mu}_m)).
 }
\end{equation*}
Now $H^0(s,\boldsymbol{\mu}_m)$ is an isomorphism while the upper right corner is trivial.
It follows that the vertical left map is trivial, so Thm.\ \ref{thm: weak section conjecture for affine brauer-severi varieties} follows from Rem.\ \ref{rem: amitsur theorem} and Lem.\ \ref{lem: brauer class and hochschild serre}.
\Endproof

\section{An analogue for the weak section conjecture: complete case.}

We are interested in an analogue of the weak section conjecture for Brauer-Severi varieties:
We want to know under what circumstances (\dag) turns out to be true for $X$ a Brauer-Severi variety over $k$.

\begin{rem}\label{rem: real case}
For any geometrically unibranch variety $Y$ over $\mathbb{R}$, such an analogue of the weak section conjecture is true without any additional assumptions:
A homotopy rational point of $Y$ gives a retraction to the canonical map $H^\bullet(\mathbb{R},\Lambda) \rightarrow H^\bullet(Y,\Lambda)$.
In particular, $Y$ has infinite cohomological dimension at $2$ and thus admits an $\mathbb{R}$-rational point by \cite{Cox79} Thm.\ 2.1.
\end{rem}

Rem.\ \ref{rem: real case} can easily be extended to any field $k$, algebraic over $\mathbb{Q}$ and henselian with respect to an archimedean valuation:

\begin{prop}\label{prop: archimedean analogue for the weak section conjecture}
Let $k/\mathbb{Q}$ be algebraic s.t.\ $k$ is henselian with respect to an archimedean valuation.
Let $Y$ be a geometrically unibranch proper $k$-variety admitting a homotopy rational point.
Then $Y$ admits a $k$-rational point.
\end{prop}

\proof
Let $v$ be the henselian archimedean valuation of $k$.
We may assume that the completion $k_v$ of $k$ is $\mathbb{R}$.
By \cite{Prestel84} Cor.\ 5.2, the ordered field $k$ is an elementary sub-structure of the ordered field $k_v$.
Thus, $X$ admits a $k$ point if and only if $X\otimes_kk_v$ admits a $k_v$-point, i.e., it suffices to show that $X\otimes_k k_v$ admits a real point (if $Y$ is a Brauer-Severi variety, this follows from ${\rm Br}(k) = {\rm Br}(k_v)$), i.e., that $X\otimes_k k_v$ has infinite cohomological dimension at $2$ (as in Rem.\ \ref{rem: real case}, see \cite{Cox79} Thm.\ 2.1).
Arguing as in Rem.\ \ref{rem: real case}, we see that $X$ has infinite cohomological dimension at $2$.
By \cite{SGA4.5} Arcata V Cor.\ 3.3, $\mathbb{R}\Gamma(X\otimes_kk^{\rm s}, \Lambda)$ and $\mathbb{R}\Gamma(X\otimes_kk_v^{\rm s}, \Lambda)$ are isomorphic in the derived category $\mathcal{D}(\underline{\rm Ab})$.
By a classical result of E.\ Artin, $k^{\rm s}/k$ is finite of degree $2$, generated by a primitive fourth root of unity (see \cite{Artin24}).
Hence, $\Gamma_k$ and $\Gamma_{k_v}$ are canonically isomorphic and $\mathbb{R}\Gamma(X\otimes_kk^{\rm s}, \Lambda)$ and $\mathbb{R}\Gamma(X\otimes_kk_v^{\rm s}, \Lambda)$ are isomorphic even in the derived category $\mathcal{D}(\underline{\rm Mod}_{\Gamma_{k_{(v)}}})$.
By taking $\Gamma_k = \Gamma_{k_v}$-hypercohomology we see that $X$ and $X\otimes_kk_v$ have isomorphic cohomology.
In particular, $X\otimes_k k_v$ has infinite cohomological dimension at $2$, just as claimed.
\Endproof

Unfortunately, Prop.\ \ref{prop: archimedean analogue for the weak section conjecture} is no longer true for $X$ a Brauer-Severi variety over a henselian non-archimedean field $k$, see Prop.\ \ref{prop: counterexample}, below.
The obstruction for an analogue of the weak section conjecture, is exactly the image of the Chern-class of a non-trivial line bundle under the homotopy rational point:

\begin{thm}\label{thm: weak section conjecture for brauer severi varieties}
Let $k$ be a field of characteristic $0$.
Let $X$ be a Brauer-Severi variety over $k$ of period $d$ with $[\mathcal{L}_X]$ a generator of ${\rm Pic}(X)$ and $s \in [B\Gamma_k,X]_{\hat{\mathcal{S}}\downarrow B\Gamma_k}$ a homotopy rational point.
\begin{enumerate}
 \item In ${\rm Br}(k)[d] = H^2(B\Gamma_k,\boldsymbol{\mu}_d)$, the relative Brauer group ${\rm Br}(X/k)$ is generated by the class $s^*\hat{c}_1[\mathcal{L}_X]$.
 \item In particular, $X$ admits no $k$-rational point, if and only if $s^*$ does not trivialize $\hat{c}_1[\mathcal{L}_X]$ modulo $d$, if and only if $s^*\hat{c}_1[\mathcal{L}_X]$ is not divisible by $d$ in $H^2(B\Gamma_k,\hat{\mathbb{Z}}(1))$.
 \item If $k$ has cohomological dimension $\leq 2$ (e.g.\ a $p$-adic or totally imaginary number field), $s$ is induced by a $k$-rational point, if and only if $s^*\hat{c}_1[\mathcal{L}_X]$ is trivial in $H^2(B\Gamma_k,\hat{\mathbb{Z}}(1))$.
\end{enumerate}
\end{thm}

\proof
First, note that (ii) is a direct consequence of (i):
Indeed, by Amitsur's Theorem, ${\rm Br}(X/k)$ is generated by the Brauer-class $[X]$ of $X$ (see e.g.\ \cite{GilleSzamuely06} Thm.\ 5.4.1).
\\
\noindent \textbf{(i):}
Under the isomorphisms ${\rm Pic}(X^{(s)}) \cong \mathbb{Z}$, the canonical projection $p:X^{\rm s}\rightarrow X$ corresponds to the multiplication-by-$d$ map.
Thus, the vertical left arrow in the commutative diagram with exact rows
\begin{equation*}
 \xymatrix{
  0 \ar[r] &
  {\rm Pic}(X) \otimes \mathbb{Z}/d\mathbb{Z} \ar[r]^-{\hat{c}_1} \ar[d]^-{\cdot d} &
  H^2(X,\boldsymbol{\mu}_d) \ar[d]^-{p^*} \ar[r] &
  {\rm Br}(X)[d] \ar[r] \ar[d] &
  0  
 \\
  0 \ar[r] &
  {\rm Pic}(X^{\rm s}) \otimes \mathbb{Z}/d\mathbb{Z} \ar[r]^-{\hat{c}_1} &
  H^2(X^{\rm s},\boldsymbol{\mu}_d) \ar[r] &
  {\rm Br}(X^{\rm s})[d] \ar[r] &
  0 
 }
\end{equation*}
is trivial.
By Lem.\ \ref{lem: second cohomology of Brauer-Severi varieties}, $H^2(X,\boldsymbol{\mu}_d)$ decomposes into $H^2(B\Gamma_k,\boldsymbol{\mu}_d)\oplus H^2(X^{\rm s},\boldsymbol{\mu}_d)$, where the projections onto the factors correspond to $s^*$ and $p^*$, so the upper kernel factors through the canonical map $H^2(B\Gamma_k,\boldsymbol{\mu}_d) \rightarrow H^2(X,\boldsymbol{\mu}_d)$:
\begin{equation*}
 \xymatrix{
  {\rm Pic}(X) \otimes \mathbb{Z}/d\mathbb{Z} \ar@{^(->}[r]^-{\hat{c}_1} \ar@{^(->}[dr] &
  H^2(X,\boldsymbol{\mu}_d) \ar[r] &
  {\rm Br}(X)\phantom{.}
 \\
  &
  H^2(B\Gamma_k,\boldsymbol{\mu}_d) \ar[u] \ar@{^(->}[r] &
  {\rm Br}(k). \ar[u]
 }
\end{equation*}
The upper row is trivial, so $\hat{c}_1$ factors even over ${\rm Br}(X/k)$.
Since the relative Brauer group is generated by the Brauer class of $X$ and $d$ is its order, (i) follows by comparing orders.
\\
\noindent \textbf{(iii):}
By Hilbert $90$, any rational point trivializes $\hat{c}_1$.
On the other hand, to show that $s^*\hat{c}_1[\mathcal{L}_X] = 0$ in $H^2(B\Gamma_k,\hat{\mathbb{Z}}(1))$ implies the rationality of $s$, first observe that $X$ splits by (ii).
It follows that up to a sign, the class $[\mathcal{L}_X]$ equals $[\mathcal{O}(1)]$ and hence $\hat{c}_1[\mathcal{L}_X]$ equals $\alpha_x$ for $x$ any $k$-rational point of $X$.
Unravelling definition, we see that $s^*\alpha_x$ is just the image of $s$ under the injection of sets of homotopy rational points studied in Prop.\ \ref{prop: injection of horat points} and $s^*\alpha_x = 0$ implies that $s$ is homotopy equivalent to $x$.
\Endproof

Let $k^\prime/k$ be a (not necessarily finite) algebraic extension.
Under the Quillen equivalence \cite{Quick10} Cor. 2.11, the base change of $k$-varieties along $k^\prime/k$ corresponds to the restriction ${\rm res}_{\Gamma_{k^\prime}}^{\Gamma_k}$ and the maps on cohomology induced by the canonical projections correspond to the maps
\begin{equation*}
 {\rm res}_{\Gamma_{k^\prime}}^{\Gamma_k}:
 \xymatrix{
  [-,K(\Lambda,q)]_{\hat{\mathcal{S}}_{\Gamma_k}} \ar[r] &
  [-\otimes_kk^\prime,K({\rm res}_{\Gamma_{k^\prime}}^{\Gamma_k}\Lambda,q)]_{\hat{\mathcal{S}}_{\Gamma_{k^\prime}}}
 }.
\end{equation*}
For $s$ a homotopy rational point of $X/k$, denote by $s_{k^\prime}$ the restricted homotopy rational point of $X\otimes_kk^\prime / k^\prime$.
\\
We are particularly interested in the following special case:
Let $k$ be a number field and let $v$ be a non-archimedean (or archimedean) valuation on $k$ together with an extension $\bar{v}$ to $k^{\rm s}$.
Denote by $k_{\bar{v}}^h$ the corresponding henselization and by $X_{\bar{v}}$ the base-change of $X$ a (Brauer-Severi) variety over $k$ to $k_{\bar{v}}^h$, i.e., Prop.\ \ref{prop: injection of horat points} and Thm.\ \ref{thm: weak section conjecture for brauer severi varieties} (iii) apply for $X_{\bar{v}} / k_{\bar{v}}^h$ in the non-archimedean case.
If $k$ admits no real archimedean place at all, Prop.\ \ref{prop: injection of horat points} and Thm.\ \ref{thm: weak section conjecture for brauer severi varieties} (iii) apply even for $X/k$.
For $s$ a homotopy rational point of $X/k$, denote by $s_{\bar{v}}$ the restricted homotopy rational point of $X_{\bar{v}} / k_{\bar{v}}^h$.
Let us draw an immediate consequence of Prop.\ \ref{prop: archimedean analogue for the weak section conjecture} and Thm.\ \ref{thm: weak section conjecture for brauer severi varieties}:

\begin{cor}\label{cor: eventually rational equals rational}
Let $k^\prime/k$ be a finite extension of a $p$-adic or totally imaginary number field $k$ and $X$ a Brauer-Severi variety over $k$ admitting a homotopy rational point $s$.
Then $s$ is rational if and only if the base extension $s_{k^\prime}$ is rational.
\end{cor}

\proof
This follows directly from Thm.\ \ref{thm: weak section conjecture for brauer severi varieties} (iii) using that the vertical right map in the commutative diagram
\begin{equation*}
 \xymatrix{
  {\rm Pic}(X\otimes_kk^\prime) \ar[r]^-{\hat{c}_1} &
  H^2(X\otimes_kk^\prime,\hat{\mathbb{Z}}(1)) \ar[r]^-{s_{k^\prime}^*} &
  H^2(\Gamma_{k^\prime}, \hat{\mathbb{Z}}(1)) 
 \\
  {\rm Pic}(X) \ar[r]^-{\hat{c}_1} \ar[u]_-{\cdot a} &
  H^2(X,\hat{\mathbb{Z}}(1)) \ar[r]^-{s^*} \ar[u] &
  H^2(\Gamma, \hat{\mathbb{Z}}(1)) \ar[u]_-{{\rm res}_{\Gamma_k}^{\Gamma_{k^\prime}}}
 }
\end{equation*}
(where $a$ is a suitable divisor of the period of $X$) is the induced map on Tate modules of the respective Brauer groups, i.e., injective (use \cite{NeukirchSchmidtWingberg} Cor.\ 7.1.4 in the local- and the Brauer-Hasse-Noether Theorem in the global case).
\Endproof

Finally, let us also mention the following two ``local-to-global''-type results:

\begin{cor}\label{cor: weak Brauer-Hasse-Noether for homotopy fixed points}
Let $X$ be a Brauer-Severi variety over $k$ a totally imaginary number field admitting a homotopy rational point.
Then the canonical map
\begin{equation*}
 \prod_v {\rm res}_{\Gamma_{k_{\bar{v}}^h}}^{\Gamma_k} :
 \xymatrix{
  [B\Gamma_k,X]_{\hat{\mathcal{S}}\downarrow B\Gamma_k} \ar[r] &
  \prod_v [B\Gamma_{k_{\bar{v}}^h},X_{\bar{v}}]_{\hat{\mathcal{S}}\downarrow B\Gamma_{k_{\bar{v}}^h}}
 }
\end{equation*}
of sets of homotopy rational points is injective.
In particular, a homotopy rational point $s$ is rational if and only if the restriction $s_{\bar{v}}$ is rational for all finite places $v$ of $k$. 
\end{cor}

\proof
Observe that $\alpha_s$ restricts to $\alpha_{s_{\bar{v}}}$, so we get the commutative diagram
\begin{equation*}
 \xymatrix{
  [B\Gamma_k,X]_{\hat{\mathcal{S}}\downarrow B\Gamma_k} \ar[rr]^-{\alpha_{s,*}} \ar[d]^-{\prod_v{\rm res}_{\Gamma_{k_{\bar{v}}^h}}^{\Gamma_k}} &&
  H^2(k,\hat{\mathbb{Z}}(1)) \ar[d]^-{\prod_v{\rm res}_{\Gamma_{k_{\bar{v}}^h}}^{\Gamma_k}} \phantom{.}
 \\
  \prod_v [B\Gamma_{k_{\bar{v}}^h},X_{\bar{v}}]_{\hat{\mathcal{S}}\downarrow B\Gamma_{k_{\bar{v}}^h}} \ar[rr]^-{\prod_v \alpha_{s_{\bar{v}},*}} &&
  \prod_v H^2(k_\nu^h,\hat{\mathbb{Z}}(1)).
 }
\end{equation*}
Thus, the claim follows from Prop.\ \ref{prop: injection of horat points} and the Brauer-Hasse-Noether Theorem.
\Endproof

\subsection*{A counterexample.}
Unfortunately, the assumption on the vanishing resp.\ divisibility of $s^*\hat{c}_1[\mathcal{L}_X]$ in Thm.\ \ref{thm: weak section conjecture for brauer severi varieties} may fail in general, i.e., there is no analogue of Prop.\ \ref{prop: archimedean analogue for the weak section conjecture} for Brauer-Severi varieties over arbitrary fields:

\begin{prop}\label{prop: counterexample}
Let $k$ be a $p$-adic field and $[A]$ a non-trivial Brauer class in the $p$-primary-torsion part ${\rm Br}(k)[p^\infty]$ of the Brauer group.
Then the corresponding Brauer-Severi variety $X_A$ admits a non-rational homotopy rational point.
\end{prop}

The idea is, to get a non-rational section out of a genus $1$ ``counterexample'' to the local weak section conjecture, i.e., out of a section $s$ of $\pi_1(C/k)$ for $C$ a genus $1$ curve without any rational points.
For the rest of the section, let $k$ be a $p$-adic field.

\begin{rem}\label{rem: relative Brauer group and maps to Brauer-Severi varieties}
Let $C$ be a connected smooth projective curve over $k$.
Further, let $[A]$ be a class in the relative Brauer group ${\rm Br}(C/k)$.
The central simple algebra $A$ splits over the function field $k(C)$ of $C$, i.e.\ $X_A$ admits a $k(C)$-rational point.
This $k(C)$-rational point extends to a non-constant rational map $C \dashrightarrow X_A$ which is even regular, since $X_A$ is proper over $k$.
Say $C$ has genus $\geq 1$, i.e., is a $K(\pi,1)$.
Then any section $s$ of $\pi_1(C/k)$ induces a homotopy rational point $f_*s$ of $X_A$ over $k$.
\end{rem}

Thus, we have to restrict our search to smooth projective curves of genus $1$ admitting no rational points, i.e., to non-split torsors under elliptic curves:

\begin{rem}\label{rem: genus 1 curves with trivial rel Brauer group are elliptic}
Let $C$ be a torsor under an elliptic curve $E$ over $k$.
Then $C$ splits if and only if the relative Brauer group ${\rm Br}(C/k)$ is trivial:
As a torsor under an elliptic curve, the canonical map $C \rightarrow {\rm Alb}^1(C)$
into the Albanese-torsor of $C$ is an isomorphism (cf.\ \cite{DescenteVI} Thm.\ 3.3).
As $C$ is a curve, ${\rm Alb}^1(C)$ is just $\underline{\rm Pic}_C^1$ and the above isomorphism is an isomorphism of $E=\underline{\rm Pic}_C^0$-torsors.
Thus, $C$ splits if and only if it has index $1$.
Finally, by \cite{Lichtenbaum69} Thm.\ 3, ${\rm Br}(C/k)$ is cyclic of order the index of $C$ ($k$ is a $p$-adic field!).
\end{rem}

\noindent \textbf{Proof of Prop.\ \ref{prop: counterexample}:}
For any given elliptic curve $E$ over $k$ there is a non-split $E$-torsor $C$ with split $\pi_1(C/k)$ by \cite{Stix12} Prop.\ 183. 
Combining Rem.\ \ref{rem: relative Brauer group and maps to Brauer-Severi varieties} and \ref{rem: genus 1 curves with trivial rel Brauer group are elliptic}, from each $[A] \in {\rm Br}(C/k)$ we get Brauer-Severi varieties $X_A$ admitting non-rational homotopy rational points.
By the genus $1$ case of \cite{Stix10} Thm.\ 2, the index of $C$ is a power of $p$, so $[A]$ is in fact $p$-primary-torsion in ${\rm Br}(k)$. 
\\
It remains to show that each $p$-primary-torsion class of ${\rm Br}(k)$ appears in the relative Brauer group for such a ``bad'' genus one curve $C$:
We start with any such non-split $E$-torsor $C$ with section $s$ of $\pi_1(C/k)$.
Let
\begin{equation*}
 \xymatrix{
  &
  C^\prime \ar[d]^-h
 \\
  B\Gamma_k \ar[ur]^-{s^\prime} \ar[r]^-s &
  C
 }
\end{equation*}
be a neighbourhood of the section $s$, i.e.\ $h$ finite \'{e}tale and $s^\prime$ a section of $\pi_1(C^\prime/k)$, compatible with $s$ via $h$.
By the Riemann-Hurwitz formula, $C^\prime$ is still a genus $1$ curve, i.e.\ a torsor under an elliptic curve with split $\pi_1(C^\prime/k)$.
As ${\rm Br}(C/k) \leq {\rm Br}(C^\prime/k)$, the colimit $\colim_{(C^\prime,s^\prime)} {\rm Br}(C^\prime/k)$ over all neighbourhoods $(C^\prime,s^\prime)$
is unbounded in ${\rm Br}(k)[p^\infty]$ by \cite{Stix12} Prop.\ 122.
Since ${\rm Br}(k)[p^\infty]$ is $\mathbb{Q}_p / \mathbb{Z}_p$, this colimit is already the whole $p$-primary part of ${\rm Br}(k)$, which finishes the proof.
\Endproof

\begin{rem}\label{rem: projective spaces have more then just one ho rat point}
Using Cor.\ \ref{cor: eventually rational equals rational} we get that the homotopy rational points of $X_A$ given by Prop.\ \ref{prop: counterexample} will never become rational after a finite extension $k^\prime/k$.
But $X_A \otimes_kk^\prime$ splits for sufficiently large $k^\prime/k$, i.e., $X_A \otimes_kk^\prime$ admits at least two different homotopy rational points: the rational and at least one non-rational one.
\end{rem}

\begin{rem}\label{rem: conjecture}
Actually, we \textit{hope} that Prop.\ \ref{prop: counterexample} is just one half of the following beautiful statement:
\textit{A Brauer-Severi variety $X$ over a $p$-adic field $k$ should admit a homotopy rational point over $k$ if and only if its period is a power of $p$.}
\end{rem}

\begin{rem}\label{rem: conjecture, known cases, ideas}
The statement of Rem.\ \ref{rem: conjecture} is true at least for homotopy rational points build out of sections of smooth projective curves of positive genus following the above recipe.
This is just \cite{Stix10} Thm.\ 2.
Unfortunately, we were not able to prove the general statement, yet.
\end{rem}

\bibliographystyle{alpha}



\bigskip{\footnotesize%
  \textsc{Johannes Schmidt, Mathematisches  Institut, Universit\"{a}t Heidelberg, Im Neuenheimer Feld 288, 69120 Heidelberg, Germany} \par  
  \textit{E-mail address}: \texttt{jschmidt@mathi.uni-heidelberg.de} \par
}

\end{document}